\documentclass{amsproc}

\usepackage[swedish,english]{babel}
\usepackage[utf8]{inputenc}
\usepackage[T1]{fontenc}
\usepackage{amsmath,amssymb,amsfonts,amsbsy}
\usepackage{url}
\usepackage[dvipsnames,usenames]{color}
\usepackage{array}
\usepackage[pdftex]{graphicx}
\usepackage{tikz}
\usepackage{accents}
\usepackage{tabu}

\newcolumntype{L}{>{\displaystyle}l}
\newcolumntype{C}{>{\displaystyle}c}
\newcolumntype{R}{>{\displaystyle}r}

\newcommand{\R}{\mathbb R}
\newcommand{\N}{\mathbb N}
\newcommand{\Z}{\mathbb Z}
\newcommand{\C}{\mathbb C}

\newcommand{\D}{\mathbb D}
\renewcommand{\Im}{\mathrm{Im}}
\renewcommand{\Re}{\mathrm{Re}}
\def\E{{\mathrm{e}}}
\def\di{\partial} 
\def\til{\widetilde}

\newcommand{\pois}{\mathcal P}

\def\I{\mathfrak{i}}
\newcommand{\diff}{\mathrm{d}}
\renewcommand{\bar}{\overline}

\newcommand{\HN}{Herglotz-Nevan\-linna }
\newcommand{\NV}{Nevanlinna }

\newcommand{\ie}{\textit{i.e.}\/ } 
\newcommand{\eg}{\textit{e.g.}\/ } 
\newcommand{\cf}{\textit{cf.}\/ } 

\renewcommand{\vec}[1]{\accentset{\rightharpoonup}{#1}} 

\allowdisplaybreaks

\theoremstyle{definition} 
\newtheorem{define}{Definition}[section]
\newtheorem{example}[define]{Example}
\newtheorem{remark}[define]{Remark}

\theoremstyle{plain} 
\newtheorem{lemma}[define]{Lemma}
\newtheorem{thm}[define]{Theorem}
\newtheorem{prop}[define]{Proposition}
\newtheorem{coro}[define]{Corollary}

\numberwithin{equation}{section}

\begin{document}

\title[Product Nevanlinna measures]{Characterizations of the Lebesgue measure and product measures related to holomorphic functions having non-negative imaginary or real part}

\author{Mitja Nedic}
\address{Mitja Nedic, Department of Mathematics, Stockholm University, SE-106 91 Sto\-ckholm, Sweden, orc-id: 0000-0001-7867-5874}
\curraddr{}
\email{mitja@math.su.se}
\thanks{ 
\\ The author was supported by the Swedish Foundation for Strategic Research, grant nr. AM13-0011.}


\date{2018-11-30} 

\begin{abstract}
In this paper, we study a class of Borel measures on $\R^n$ that arises as the class of representing measures of \HN functions. In particular, we study product measures within this class where products with the Lebesgue measures play a special role. Hence, we give several characterizations of the $n$-dimensional Lebesgue measure among all such measures and characterize all product measures that appear in this class of measures. Furthermore, analogous results for the class of positive Borel measures on the unit poly-torus with vanishing mixed Fourier coefficients are also presented, and the relation between the two classes of measures with regard to the obtained results is discussed.
\end{abstract}

\maketitle

\section{Introduction}\label{sec:intro}

Holomorphic functions with a prescribed sign of their imaginary or real part are of great interest in complex analysis and appear in many areas and applications. To name but a few examples, such functions of one variable appear when considering the moment problem \cite{Akhiezer1965,Nevanlinna1922,Simon1998}, in spectral theory of Sturm-Liouville problems and perturbations \cite{Aronszajn1957,AronszajnBrown1970,Donoghue1965,KacKrein1974}, when deriving physical bounds for passive systems \cite{BernlandEtal2011} or as approximating functions in certain convex optimization problems \cite{IvanenkoETAL2019a,IvanenkoETAL2018}. On the other hand, such functions of several variables are connected \eg with operator monotone functions \cite{AglerEtal2012} or with representations of multidimesional passive systems \cite{Vladimirov1979}.

With respect to the domain of the function in question, the two most common cases are to consider the poly-upper half-plane and the unit polydisk. In the case of the unit polydisk, one may consider functions with non-negative real part, \cf \cite{KoranyiPukanszky1963} and Section \ref{sec:polydisk}, or with non-negative imaginary part as in \eg \cite{VladimirovDrozzinov1974}. However, in this paper, it is the case poly-upper half-plane that will be of primary interest. We note that more general domains are considered in \eg \cite{AizenbergDautov1976,AizenbergYuzhakov1983,Vladimirov1969}.

In the poly-upper half-plane $\C^{+n}$, \ie 
$$\C^{+n}:= (\C^+)^n = \big\{\vec{z}\in\C^n \,\big |\,\forall j=1,2,\ldots, n:  \Im[z_j]>0 \big\},$$
we consider the class of all holomorphic functions with non-negative imaginary part. This is a well-studied class of functions, appearing, as a whole or as one of its subclasses, in \eg \cite{AglerEtal2012,AglerEtal2016,LugerNedic2017,LugerNedic2018,LugerNedic2019,Savchuk2006,Vladimirov1969,Vladimirov1979}. These functions, called \emph{Herglotz-Nevanlinna functions}, \cf Definition \ref{def:HN_functions}, admit a power-full integral representation theorem, \cf Theorem \ref{thm:intRep_Nvar}, that characterizes this class of functions in terms of a real number $a$, a vector $\vec{b} \in [0,\infty)^n$ and a positive Borel measure $\mu$ on $\R^n$ satisfying two conditions. It is this class of representing measures, called \emph{Nevanlinna measures}, \cf Definition \ref{def:Nevan}, that is the main object of study in this paper.

As one of the representing parameters of \HN functions, Nevanlinna measures have been studied in \cite{LugerNedic2017,LugerNedic2018,LugerNedic2019}. While these measures are somewhat unremarkable in dimension 1, where they simply consists of all positive Borel measure $\mu$ on $\R$ such that $\int_\R(1+t^2)^{-1}\diff\mu(t) < \infty$, in higher dimension, the class becomes much more obscure. For example, it is known that such measures cannot be finite unless trivial \cite[Prop. 4.3]{LugerNedic2017}, they cannot have point masses \cite[Prop. 4.4]{LugerNedic2017}, their restriction to coordinate parallel hyperplanes is very particular \cite[Thm. 3.4]{LugerNedic2018} and their support must obey certain geometric restrictions \cite[Thms. 3.10, 3.16 and 3.24]{LugerNedic2018}.

In this paper, our main focus will be on the subclass of product Nevanlinna measures, \ie Nevanlinna measures that are product measures. In particular, products with the Lebesgue measures are of great interest as they represent functions that do not depend on all of their variables and can be interpolated by convex combinations, \cf Section \ref{sec:product_measures} and \cite{LugerNedic2019,Nedic2017}. Therefore, our objective is twofold. Firstly, we would like to give conditions that characterize the $n$-dimensional Lebesgue measure among all Nevanlinna measures on $\R^n$. This is presented in Theorem \ref{thm:Lebesgue_complete}, see also Propositions \ref{prop:Lebesgue_v2} and \ref{prop:Lebesgue}. Due to their particular interest, the low-dimension cases are also listed separately as corollaries of the main theorem. Secondly, we would like to identify all Nevanlinna measures that are product measures, which is presented in Theorem \ref{thm:product_Nevan} and Corollary \ref{coro:mixed_product}.

Due to their close connection with Nevanlinna measures, measures on the unit poly-torus, or, equivalently, on $[0,2\pi)^n$, that have vanishing mixed Fourier coefficients are also discussed. Analogous results are presented, \ie several characterizations of the Lebesgue measure among all such measures are given in Theorem \ref{thm:Lebesgue_complete_polydisk}, while product measures are identified in Proposition \ref{prop:product_polydisk}. 

The structure of the paper is as follows. After the introduction in Section \ref{sec:intro} we recall the necessary prerequisites regarding \HN functions and Nevanlinna measures in Section \ref{sec:HN_functions}. Section \ref{sec:Lebesgue} is devoted to the results concerning characterizations of the Lebesgue measure. Product measures are, afterwards, discussed in Section \ref{sec:product_measures}. Finally, Section \ref{sec:polydisk} discusses measures on $[0,2\pi)^n$ with vanishing mixed Fourier coefficients.

\section{Herglotz-Nevanlinna functions and Nevanlinna measures}\label{sec:HN_functions}

The class of \HN functions on the poly-upper half-plane is formally defined as follows \cite{LugerNedic2017,LugerNedic2019,Vladimirov1979}.

\begin{define}\label{def:HN_functions}
A function $q\colon \C^{+n} \to \C$ is called a \emph{\HN function} if it is holomorphic and has non-negative imaginary part.
\end{define}

The primary tool in the study of \HN functions is the following integral representation theorem, \cf \cite[Thm. 4.1 and Thm. 5.1]{LugerNedic2019}.

\begin{thm}\label{thm:intRep_Nvar}
A function $q\colon \C^{+n} \to \C$ is a \HN function if and only if $q$ can be written, for every $\vec{z} \in \C^{+n}$, as
\begin{equation}\label{eq:intRep_Nvar}
q(\vec{z}) = a + \sum_{\ell=1}^nb_\ell\,z_\ell + \frac{1}{\pi^n}\int_{\R^n}K_n(\vec{z},\vec{t})\diff\mu(\vec{t}),
\end{equation}
where $a \in \R$, $\vec{b} \in [0,\infty)^n$, the kernel $K_n$ is defined for $\vec{z} \in \C^{+n}$ and $\vec{t} \in \R^n$ as
$$K_n(\vec{z},\vec{t}) := \I\left(\frac{2}{(2\I)^n}\prod_{\ell=1}^n\left(\frac{1}{t_\ell-z_\ell}-\frac{1}{t_\ell+\I}\right)-\frac{1}{(2\I)^n}\prod_{\ell=1}^n\left(\frac{1}{t_\ell-\I}-\frac{1}{t_\ell+\I}\right)\right)$$
and $\mu$ is a positive Borel measure on $\R^n$ satisfying the growth condition
\begin{equation}
\label{eq:measure_growth}
\int_{\R^n}\prod_{\ell=1}^n\frac{1}{1+t_\ell^2}\diff\mu(\vec{t}) < \infty
\end{equation}
and the \NV condition
\begin{equation}
\label{eq:measure_Nevan}
\int_{\R^n}\frac{1}{(t_{\ell_1} - z_{\ell_1})^2(t_{\ell_2} - \bar{z_{\ell_1}})^2}\prod_{\substack{j=1 \\ j \neq \ell_1,\ell_2}}^n\left(\frac{1}{t_j - z_j} - \frac{1}{t_j - \bar{z}_j}\right)\diff\mu(\vec{t}) = 0
\end{equation}
for all $\vec{z} \in \C^{+n}$ and all indices $\ell_1,\ell_2 \in \{1,2,\ldots,n\}$ with $\ell_1 < \ell_2$. Furthermore, for a given function $q$, the triple of representing parameters $(a,\vec{b},\mu)$ is unique.
\end{thm}

This theorem gives rise to the following class of measures on $\R^n$, \cf \cite[Def. 3.1]{LugerNedic2018}.

\begin{define}\label{def:Nevan}
A positive Borel measure $\mu$ on $\R^n$ is called a \emph{Nevanlinna measure} if it is the representing measure of some \HN function, \ie if it satisfies the growth condition \eqref{eq:measure_growth} and the Nevanlinna condition \eqref{eq:measure_Nevan}.
\end{define}

One of the most important Nevanlinna measures is the Lebesgue measure, which is connected to the following \HN functions, \cf \cite[Ex. 3.5]{LugerNedic2019}.

\begin{example}\label{ex:Lebesgue}
Let $n \in \N$ and consider the function $q\colon \C^{+n} \to \C$, given by
$$q(\vec{z}) := \I.$$
This function is represented by the data $(0,\vec{0},\lambda_{\R^n})$ in the sense of representation \eqref{eq:intRep_Nvar}, where $\lambda_{\R^n}$ denotes the standard Lebesgue measure on $\R^n$. More generally, the function
$$q(\vec{z}) := \eta$$
for some number $\eta \in \C^+ \cup \R$ is represented by the data $(\Re[\eta],\vec{0},\Im[\eta]\lambda_{\R^n})$.\hfill$\lozenge$
\end{example}

It is known that the \NV measures may be characterized via different conditions, two of which are expressed in terms of the factors that make up the kernel $K_n$. To that end, we introduce, for $k \in \{-1,0,1\}$ and $j \in \{1,2,\ldots,n\}$, the factors $N_{k,j}$ similarly to \cite[pg. 1193]{LugerNedic2019}. Given ambient vectors $\vec{z} \in \C^{+n}$ and $\vec{t} \in \R^n$, we define
\begin{eqnarray*}
N_{-1,j} & := & \frac{1}{2\,\I}\left(\frac{1}{t_j - z_j} - \frac{1}{t_j - \I}\right), \\
N_{0,j} & := & \frac{1}{2\,\I}\left(\frac{1}{t_j - \I} - \frac{1}{t_j + \I}\right), \\
N_{1,j} & := & \frac{1}{2\,\I}\left(\frac{1}{t_j + \I} - \frac{1}{t_j - \bar{z}_j}\right),
\end{eqnarray*}
leading to the following equivalence \cite[Thm. 5.1]{LugerNedic2019}.

\begin{thm}
    Let $n \geq 2$ and let $\mu$ be a positive Borel measure on $\R^n$ satisfying the growth condition \eqref{eq:measure_growth}. Then the following statements are equivalent.
    \begin{itemize}
    \item[(a)]{It holds that
    \begin{equation}
        \label{eq:Nevan_original}
        \sum_{\substack{\vec{\rho} \in \{-1,0,1\}^n \\ -1\in\vec{\rho}\,\wedge\,1\in\vec{\rho}}}\int_{\R^n}N_{\rho_1,1}N_{\rho_2,2}\ldots N_{\rho_n,n}\diff\mu(\vec{t}) = 0
    \end{equation}
    for all $\vec{z} \in \C^{+n}$.}\vspace{1mm}
    \item[(b)]{It holds that
    \begin{equation}
        \label{eq:Nevan_separate}
        \int_{\R ^n}N_{\rho_1,1}N_{\rho_2,2}\ldots N_{\rho_n,n}\diff\mu(\vec{t}) = 0
    \end{equation}
    for all $\vec{z} \in \C^{+n}$ and for every vector $\vec{\rho} \in \{-1,0,1\}^n$ with at least one entry equal to $1$ and at least one entry equal to $-1$.}\vspace{1mm}
    \item[(c)]{The measure $\mu$ satisfies the Nevanlinna condition \eqref{eq:measure_Nevan}.}\vspace{1mm}
    \item[(d)]{It holds that
    \begin{equation}
        \label{eq:Nevan_discrete}
        \int_{\R ^n}\left(\frac{t_1-\I}{t_1+\I}\right)^{m_1}\ldots\:\left(\frac{t_n-\I}{t_n+\I}\right)^{m_n}\prod_{\ell=1}^n\frac{1}{1+t_\ell^2}\diff\mu(\vec{t}) = 0
    \end{equation}
    for all multi-indices $\vec{m} \in \Z^n$ with at least one positive entry and at least one negative entry.}
    \end{itemize}
\end{thm}

\section{Characterizations of the Lebesgue measure}\label{sec:Lebesgue}

Our first main theorem provides four characterizations of the Lebesgue measure among all Nevanlinna measures.

\begin{thm}\label{thm:Lebesgue_complete}
Let $n \in \N$ and let $\mu$ be a Nevanlinna measure in $\R^n$. Then, $\mu = c \lambda_{\R^n}$ for some constant $c \geq 0$ if and only if the measure $\mu$ satisfies one (and thus all) of the following conditions:
\begin{itemize}
    \item[(a)]{It holds that
    \begin{equation}
        \label{eq:Lebesgue_Nevan_original}
        \sum_{\substack{\vec{\rho} \in \{-1,0,1\}^n \\ 1\not\in\vec{\rho}\,\wedge\,\vec{\rho} \neq \vec{0}}}\int_{\R^n}N_{\rho_1,1}N_{\rho_2,2}\ldots N_{\rho_n,n}\diff\mu(\vec{t}) = 0
    \end{equation}
    for all $\vec{z} \in \C^{+n}$.    
    }
    \item[(b)]{It holds that
    \begin{equation}
        \label{eq:Lebesgue_Nevan_separate}
        \int_{\R^n}N_{\rho_1,1}N_{\rho_2,2}\ldots N_{\rho_n,n}\diff\mu(\vec{t}) = 0
    \end{equation}
    for all $\vec{z} \in \C^{+n}$ and and for every vector $\vec{\rho} \in \{-1,0,1\}^n$ that is not the zero-vector and has no entries equal to $1$.
    }
    \item[(c)]{It holds, for every index $j \in \{1,2,\ldots,n\}$, that
    \begin{equation}
        \label{eq:Lebesgue}
        \int_{\R^n}\frac{1}{(t_j - z_j)^2}\prod_{\substack{\ell = 1 \\ \ell \neq j}}^n\left(\frac{1}{t_\ell - z_\ell} - \frac{1}{t_\ell - \bar{z}_\ell}\right)\diff\mu(\vec{t}) = 0
    \end{equation}
    for all $\vec{z} \in \C^{+n}$.
    }
    \item[(d)]{It holds that
    $$\int_{\R^n}\left(\frac{t_1-\I}{t_1+\I}\right)^{m_1}\ldots\left(\frac{t_n-\I}{t_n+\I}\right)^{m_n}\prod_{j=1}^n\frac{1}{1+t_j^2}\diff\mu(\vec{t}) = 0$$
    for every multi-index $\vec{m} \in (\N_0)^n\setminus\{\vec{0}\}$.
    } 
\end{itemize}
\end{thm}

First, statements (a) and (b) will be established in Proposition \ref{prop:Lebesgue_v2} and Corollary \ref{coro:Lebesgue_v2}. Afterwards, statement (c) will be established in Proposition \ref{prop:Lebesgue}.  Finally, statement (d) will follow from the corresponding result on the unit polydisk in Proposition \ref{prop:Lebesgue_polydisk}.

A reformulation of Theorem \ref{thm:Lebesgue_complete} may also be stated as follows.

\begin{coro}
    Let $n \in \N$ and let $\mu$ be a positive Borel measure on $\R^n$ satisfying the growth condition \eqref{eq:measure_growth}. Then, $\mu = c\lambda_{\R^n}$ for some constant $c \geq 0$ if and only if the measure $\mu$ satisfies one (and thus all) of the following conditions:
    \begin{itemize}
    \item[(a)]{It holds that
    \begin{equation}
    \label{eq:Lebesgue_Nevan_complete}
    \sum_{\substack{\vec{\rho} \in \{-1,0,1\}^n \\ \vec{\rho} \neq \vec{0}}}\int_{\R^n}N_{\rho_1,1}N_{\rho_2,2}\ldots N_{\rho_n,n}\diff\mu(\vec{t}) = 0
    \end{equation}
    for all $\vec{z} \in \C^{+n}$.}\vspace{1mm}
    \item[(b)]{It holds that
        $$\int_{\R ^n}N_{\rho_1,1}N_{\rho_2,2}\ldots N_{\rho_n,n}\diff\mu(\vec{t}) = 0$$
    for all $\vec{z} \in \C^{+n}$ and for every vector $\vec{\rho} \in \{-1,0,1\}^n$ that is not the zero-vector.}\vspace{1mm}
    \item[(c)]{The measure $\mu$ satisfies the Nevanlinna condition \eqref{eq:measure_Nevan} and condition \eqref{eq:Lebesgue}.}\vspace{1mm}
    \item[(d)]{It holds that
    $$\int_{\R^n}\left(\frac{t_1-\I}{t_1+\I}\right)^{m_1}\ldots\left(\frac{t_n-\I}{t_n+\I}\right)^{m_n}\prod_{j=1}^n\frac{1}{1+t_j^2}\diff\mu(\vec{t}) = 0$$
    for every multi-index $\vec{m} \in \Z^n\setminus\{\vec{0}\}$.}
    \end{itemize}
\end{coro}

If one wants to characterize the standard Lebesgue measure via the above theorem, \ie the case when $c=1$, then the following normalization constraint should be added.

\begin{coro}
    Let $n \in \N$ and let $\mu$ be a Nevanlinna measure in $\R^n$. Then, $\mu = \lambda_{\R^n}$ if and only if $\mu$ satisfies one (and thus all) of the condition of Theorem \ref{thm:Lebesgue_complete} and
    $$\int_{\R^n}\frac{1}{1+t_\ell^2}\diff\mu(\vec{t}) = \pi^n.$$
\end{coro}

Before proceeding, we also state the variants of Theorem \ref{thm:Lebesgue_complete} for dimensions $n=1$ and $n=2$ as separate corollaries due to their particular interest.

\begin{coro}\label{coro:Lebesgue_dim1}
    Let $\mu$ be a positive Borel measure on $\R$ satisfying the growth condition \eqref{eq:measure_growth}. Then, $\mu = c\lambda_{\R}$ for some constant $c \geq 0$ if and only if the measure $\mu$ satisfies one (and thus all) of the following conditions:
    \begin{itemize}
    \item[(a)]{It holds that
    $$\int_\R\left(\frac{\Im[z]}{|t-z|^2}-\frac{1}{1+t^2}\right)\diff\mu(t) = 0$$
    for all $z \in \C^+$.}\vspace{1mm}
    \item[(b)]{It holds that
    $$\int_\R\left(\frac{1}{t-z} - \frac{1}{t-\I}\right)\diff\mu(t) = 0$$
    for all $z \in \C^+$.}\vspace{1mm}
    \item[(c)]{It holds that
    $$\int_\R\frac{1}{(t-z)^2}\diff\mu(t) = 0$$
    for all $z \in \C^+$.}\vspace{1mm}
    \item[(d)]{It holds that
    $$\int_{\R}\left(\frac{t-\I}{t+\I}\right)^{m}\frac{1}{1+t^2}\diff\mu(t) = 0$$
    for every index $m \in \Z\setminus\{0\}$.}
    \end{itemize}
\end{coro}

\begin{coro}\label{coro:Lebesgue_dim2}
    Let $\mu$ be a positive Borel measure on $\R^2$ satisfying the growth condition \eqref{eq:measure_growth}. Then, $\mu = c\lambda_{\R^2}$ for some number $c \geq 0$ if and only if the measure $\mu$ satisfies one (and thus all) of the following conditions:
    \begin{itemize}
    \item[(a)]{It holds that
    $$\iint_{\R^2}\left(\frac{\Im[z_1]}{|t_1-z_1|^2}\,\frac{\Im[z_2]}{|t_2-z_2|^2}-\frac{1}{1+t_1^2}\,\frac{1}{1+t_2^2}\right)\diff\mu(t_1,t_2) = 0$$
    for all $(z_1,z_2) \in \C^{+2}$.}\vspace{1mm}
    \item[(b)]{It holds that
    \begin{eqnarray*}
    \iint_{\R^2}\left(\frac{1}{t_1-z_1} - \frac{1}{t_1-\I}\right)\left(\frac{1}{t_2-\I} - \frac{1}{t_2+\I}\right)\diff\mu(t_1,t_2) & = & 0, \\
    \iint_{\R^2}\left(\frac{1}{t_1-\I} - \frac{1}{t_1+\I}\right)\left(\frac{1}{t_2-z_2} - \frac{1}{t_2-\I}\right)\diff\mu(t_1,t_2) & = & 0, \\
    \iint_{\R^2}\left(\frac{1}{t_1-z_1} - \frac{1}{t_1-\I}\right)\left(\frac{1}{t_2-z_2} - \frac{1}{t_2-\I}\right)\diff\mu(t_1,t_2) & = & 0, \\
    \iint_{\R^2}\left(\frac{1}{t_1-z_1} - \frac{1}{t_1-\I}\right)\left(\frac{1}{t_2-\bar{z}_2} - \frac{1}{t_2+\I}\right)\diff\mu(t_1,t_2) & = & 0, \\
    \end{eqnarray*}
    for all $(z_1,z_2) \in \C^{+2}$.}\vspace{1mm}
    \item[(c)]{It holds that
    it holds that
    \begin{eqnarray*}
    \iint_{\R^2}\frac{1}{(t_1-z_1)^2(t_2-\bar{z}_2)^2}\diff\mu(t_1,t_2) & = & 0, \\
    \iint_{\R^2}\frac{1}{(t_1-z_1)^2}\left(\frac{1}{t_2-z_2} - \frac{1}{t_2-\bar{z}_2}\right)\diff\mu(t_1,t_2) & = & 0, \\
    \iint_{\R^2}\left(\frac{1}{t_1-z_1} - \frac{1}{t_1-\bar{z}_1}\right)\frac{1}{(t_2-\bar{z}_2)^2}\diff\mu(t_1,t_2) & = & 0,
    \end{eqnarray*}
    for all $(z_1,z_2) \in \C^{+2}$.}\vspace{1mm}
    \item[(d)]{It holds that
    $$\iint_{\R^2}\left(\frac{t_1-\I}{t_1+\I}\right)^{m_1}\left(\frac{t_2-\I}{t_2+\I}\right)^{m_2}\prod_{j=1}^2\frac{1}{1+t_j^2}\diff\mu(t_1,t_2) = 0$$
    for every multi-index $\vec{m} \in \Z^2\setminus\{\vec{0}\}$.}
    \end{itemize}
\end{coro}

\subsection{First and second characterizations}\label{sec:Lebesgue_v2}

Nevanlinna measures have an intricate connection to the Poisson kernel $\pois_n$ of $\C^{+n}$, which, we recall, can be written using complex coordinates as
$$\pois_n(\vec{z},\vec{t}) := \prod_{\ell=1}^n\frac{\Im[z_\ell]}{|t_\ell - z_\ell|^2} = \frac{1}{(2\I)^n}\prod_{\ell=1}^n\left(\frac{1}{t_\ell-z_\ell}-\frac{1}{t_\ell-\bar{z_\ell}}\right).$$
Note that $\pois_n > 0$ for every $\vec{z} \in \C^{+n}$ and $\vec{t} \in \R^n$. 
Conditions \eqref{eq:Nevan_original} and \eqref{eq:Nevan_separate} link the imaginary part of the kernel $K_n$ to the Poisson kernel $\pois_n$, ensuring, thereby, that the imaginary part of representation \eqref{eq:intRep_Nvar} is non-negative. Indeed, due to \cite[Prop. 3.3]{LugerNedic2019}, conditions \eqref{eq:Nevan_original} and \eqref{eq:Nevan_separate} ensure that
$$\int_{\R^n}\Im[K_n(\vec{z},\vec{t})]\diff\mu(\vec{t}) = \int_{\R^n}\pois_n(\vec{z},\vec{t})\diff\mu(\vec{t}) \geq 0.$$

We may now characterize the Lebesgue measure via a condition that reflects the description of Nevanlinna measures by conditions \eqref{eq:Nevan_original} and \eqref{eq:Nevan_separate}, \ie statements (a) and (b) from Theorem \ref{thm:Lebesgue_complete}.

\begin{prop}\label{prop:Lebesgue_v2}
    Let $n \in \N$ and let $\mu$ be a Nevanlinna measure in $\R^n$. Then, $\mu = c\lambda_{\R^n}$ for some constant $c \geq 0$ if and only if it satisfies condition \eqref{eq:Lebesgue_Nevan_original}, \ie it holds that
    \begin{equation}\tag{\ref{eq:Lebesgue_Nevan_original}}
        \sum_{\substack{\vec{\rho} \in \{-1,0,1\}^n \\ 1\not\in\vec{\rho}\,\wedge\,\vec{\rho} \neq \vec{0}}}\int_{\R^n}N_{\rho_1,1}N_{\rho_2,2}\ldots N_{\rho_n,n}\diff\mu(\vec{t}) = 0
    \end{equation}
    for all $\vec{z} \in \C^{+n}$.
\end{prop}

\proof
Observe first that
$$\int_\R N_{-1,j}\diff t_j = 0$$
for every $z \in \C^+$ by standard residue calculus, yielding that condition \eqref{eq:Lebesgue_Nevan_original} is satisfied by the measure $c\lambda_\R$ as each index-vector $\vec{\rho}$ that appears in the sum in condition \eqref{eq:Lebesgue_Nevan_original} has at lest one entry equal to $-1$.

Conversely, let $\mu$ be a Nevanlinna measure that satisfies, in addition, condition \eqref{eq:Lebesgue_Nevan_original}. Let $q$ be the \HN function given by the data $(0,0,\mu)$ in the sense of Theorem \ref{thm:intRep_Nvar}. We now claim that, for this function $q$, it holds that $\Re[q] \equiv 0$.

To that end, note that the kernel $K_n$ may be written in terms of the factors $N_{k,j}$ as
$$K_n = \I\bigg(2\prod_{j=1}^n(N_{-1,j}+N_{0,j}) - \prod_{j=1}^nN_{0,j}\bigg).$$
Using the identities that $\bar{N}_{-1,j} = N_{1,j}$ and $\bar{N}_{0,j} = N_{0,j}$, we calculate that
\begin{multline*}
\Re[K_n] = \I\bigg(\prod_{j=1}^n(N_{-1,j}+N_{0,j})-\prod_{j=1}^n(N_{1,j}+N_{0,j})\bigg) \\
= \I\bigg(\sum_{\substack{\vec{\rho} \in \{-1,0,1\}^n \\ 1\not\in\vec{\rho}\,\wedge\,\vec{\rho} \neq \vec{0}}}N_{\rho_1,1}N_{\rho_2,2}\ldots N_{\rho_n,n} - \sum_{\substack{\vec{\rho} \in \{-1,0,1\}^n \\ -1\not\in\vec{\rho}\,\wedge\,\vec{\rho} \neq \vec{0}}}N_{\rho_1,1}N_{\rho_2,2}\ldots N_{\rho_n,n}\bigg).
\end{multline*}
As the measure $\mu$ satisfies condition \eqref{eq:Lebesgue_Nevan_original}, it also satisfies the same condition where all of the $N_{k,j}$-factors have been conjugated. Therefore, it holds, for any $\vec{z} \in \C^{+n}$, that
$$\Re[q(\vec{z})] = \int_{\R^n}\Re[K_n(\vec{z},\vec{t})]\diff\mu(\vec{t}) = 0.$$
However, a holomorphic function on a simply connected domain whose real part is identically zero must be equal to pure-imaginary constant function, \ie there exists a constant $c \geq 0$ such that $q(\vec{z}) = c\,\I$ for all $\vec{z} \in \C^{+n}$. By the uniqueness statement of Theorem \ref{thm:intRep_Nvar} and Example \ref{ex:Lebesgue}, the result follows.
\endproof

The following corollary is now an immediate consequence of the above proposition.

\begin{coro}\label{coro:Lebesgue_v2}
    Let $n \in \N$ and let $\mu$ be a Nevanlinna measure in $\R^n$. Then, the measure $\mu$ satisfies condition \eqref{eq:Lebesgue_Nevan_original} if and only if it satisfies condition \eqref{eq:Lebesgue_Nevan_separate}, \ie it holds that
    \begin{equation}\tag{\ref{eq:Lebesgue_Nevan_separate}}
        \int_{\R ^n}N_{\rho_1,1}N_{\rho_2,2}\ldots N_{\rho_n,n}\diff\mu(\vec{t}) = 0
    \end{equation}
    for all $\vec{z} \in \C^{+n}$ and for every vector $\vec{\rho} \in \{-1,0,1\}^n$ that is not the zero-vector and has no entries equal to $1$.
\end{coro}

If conditions \eqref{eq:Nevan_original} and \eqref{eq:Lebesgue_Nevan_original} are combined as in condition \eqref{eq:Lebesgue_Nevan_complete}, a further connection to the Poisson kernel becomes visible. Due to how the Poisson kernel may be written in terms of the $N_{k,j}$-factors, \cf \cite[Prop. 3.3]{LugerNedic2019} and its proof, it holds that
$$\sum_{\substack{\vec{\rho} \in \{-1,0,1\}^n \\ \vec{\rho} \neq \vec{0}}}N_{\rho_1,1}N_{\rho_2,2}\ldots N_{\rho_n,n} = \pois_n(\vec{z},\vec{t}) - \prod_{j=1}^nN_{0,j} = \pois_n(\vec{z},\vec{t}) - \pois_n(\I\,\vec{1},\vec{t})$$
for every $\vec{z} \in \C^{+n}$ and every $\vec{t} \in \R^n$. This was also seen explicitly in Corollaries \ref{coro:Lebesgue_dim1} and \ref{coro:Lebesgue_dim2}.

\subsection{Third characterization}\label{sec:Lebesgue_v1}

The \NV condition \eqref{eq:measure_Nevan} assures that, for any Nevanlinna measure $\mu$, the function
\begin{equation}
    \label{eq:Poisson_function}
    \vec{z} \mapsto \int_{\R^n}\pois_n(\vec{z},\vec{t})\diff\mu(\vec{t})
\end{equation}
is a pluriharmonic function on $\C^{+n}$, \cf \cite[Thm. 5.1]{LugerNedic2019}. Indeed, one may recognize the factors that make up the integrand in condition \eqref{eq:measure_Nevan} as being the same factors that appear in the definition of $\pois_n$ or their derivatives, assuring, thus, that the appropriate derivatives of the function \eqref{eq:Poisson_function} are identically zero.

Statement (c) of Theorem \ref{thm:Lebesgue_complete} may, therefore, be thought of as a characterization of the Lebesgue measure that uses the derivatives of the Poisson kernel. In order to prove this, we will require additional information, namely two results concerning the symmetric extension of a \HN function and a combinatorial lemma.

For the former, we recall that the integral representation in formula \eqref{eq:intRep_Nvar} is well-defined for any $\vec{z} \in (\C\setminus\R)^n$, which may be used to extend any Herglotz-Nevanlinna function $q$ from $\C^{+n}$ to $(\C\setminus\R)^n$. This extension is called the \emph{symmetric extension} of the function $q$ and is denoted as $q_\mathrm{sym}$. We note that that the symmetric extension of a \HN function $q$ is different from its possible analytic extension as soon as $\mu \neq 0$, \cf \cite[Prop. 6.10]{LugerNedic2019}. The symmetric extension also satisfies the following variable-dependence property, \cf \cite[Prop. 6.9]{LugerNedic2019}.

\begin{prop}\label{prop:variable_dependence}
    Let $n \geq 2$ and let $q_\mathrm{sym}$ be the symmetric extension of a Herglotz-Nevanlinna function $q$ in $n$ variables. Let $\vec{z} \in (\C\setminus\R)^n$ bu such that $z_j \in \C^-$ for some index $j \in \{1,2,\ldots,n\}$. Then, the value $q_\mathrm{sym}(\vec{z})$ does not depend on the components of $\vec{z}$ that lie in $\C^+$. In particular, for $\ell \neq j$, it holds that
    $$\frac{\di q_\mathrm{sym}}{\di z_\ell}(\vec{z}) = 0$$
    whenever $z_j \in \C^-$ and $z_\ell \in \C^+$.
\end{prop}

As the \HN function itself, its partial derivatives also admit an integral representation formula. To describe them, we use standard multi-index notation. Let $\vec{k} = (k_1,k_2,\ldots,k_n) \in \N_0^n$ be a multi-index with
$$|\vec{k}| := k_1 + k_2 + \ldots + k_n.$$
If $|\vec{k}| \geq 1$, the partial derivatives of the kernel $K_n$ are given by
\begin{multline}
    \label{eq:kernel_derivative}
    \frac{\di^{|\vec{k}|}}{\di z_1^{k_1} \di z_2^{k_2} \ldots \di z_n^{k_n}}K_n(\vec{z},\vec{t}) \\ = \frac{1}{(2\,\I)^{n-1}}\prod_{\substack{\ell = 1 \\ k_\ell = 0}}^n\left(\frac{1}{t_\ell - z_\ell} - \frac{1}{t_\ell + \I}\right)\prod_{\substack{\ell = 1 \\ k_\ell > 0}}^n\frac{1}{(t_\ell-z_\ell)^{k_\ell+1}},
\end{multline}
leading to the following representation formula for the partial derivatives of (the symmetric extension of) a \HN function.

\begin{prop}
\label{prop:HN_derivatives}
Let $n \geq 1$, let $q_\mathrm{sym}$ be the symmetric extension of a \HN function $q$ represented by the data $(a,\vec{b},\mu)$ and let $\vec{k} \in \N_0^n$ be a multi-index with $|\vec{k}| \geq 1$. Then, it holds, for $\vec{z} \in (\C\setminus\R)^n$, that
\begin{multline}
    \label{eq:HN_derivatives}
    \frac{\di^{|\vec{k}|}}{\di z_1^{k_1} \di z_2^{k_2} \ldots \di z_n^{k_n}}q_\mathrm{sym}(\vec{z}) = \sum_{\ell=1}^n\left\{\begin{array}{rcl}
    b_\ell & ; & |\vec{k}| = 1 \wedge k_\ell = 1, \\
    0 & ; & \text{otherwise},
    \end{array}\right. \\
    + \frac{1}{\pi^{n}}\int_{\R^n}\frac{\di^{|\vec{k}|}}{\di z_1^{k_1} \di z_2^{k_2} \ldots \di z_n^{k_n}}K_n(\vec{z},\vec{t})\diff\mu(\vec{t}).
\end{multline}
\end{prop}

\proof
In order to prove the above proposition, it suffices to prove that we may interchange taking the derivative and integration with respect to the measure $\mu$. To that end, we note that the kernel $K_n$ is uniformly bounded on any compact set $U \subseteq (\C\setminus\R)^n$. Indeed, the kernel $K_n$ may be rewritten as
$$K_n(\vec{z},\vec{t}) = \frac{\I^{3n+1}\prod_{j=1}^n(t_j - \I)(z_j + \I) - 2^{n-1}\I\prod_{j = 1}^n(t_j-z_j)}{2^{n-1}\prod_{j = 1}^n(t_j-z_j)}\:\prod_{j = 1}^n\frac{1}{1+t_j^2}$$
and for every compact set $U \subseteq (\C\setminus\R)^n$ one can find a constant $C$ (depending on $U$), such that
$$\left|\frac{\I^{3n+1}\prod_{j=1}^n(t_j - \I)(z_j + \I) - 2^{n-1}\I\prod_{j = 1}^n(t_j-z_j)}{2^{n-1}\prod_{j = 1}^n(t_j-z_j)}\right| \leq C$$
for every $\vec{z} \in U$ and $\vec{t} \in \R^n$. This finishes the proof.
\endproof

Finally, we will need the following combinatorial lemma.

\begin{lemma}\label{lem:combinatorial}
    Let $\rho \in \{0,1\}$ and $\psi_\rho\colon \C \to \C$ be the function
    $$\psi_\rho(z) := \left\{\begin{array}{rcl}
    z & ; & \rho = 0, \\
    \bar{z} & ; & \rho = 1.
    \end{array}\right.$$
    Let $n \in \N$ and let $\vec{\xi},\vec{\eta} \in \C^n$. Then, it holds that
    \begin{equation}
        \label{eq:combinatorial}
        \sum_{\vec{\rho} \in \{0,1\}^n}(-1)^{|\vec{\rho}|}\prod_{j = 1}^n(\psi_{\rho_j}(\xi_j) - \eta_j) = \prod_{j = 1}^n(\xi_j - \bar{\xi}_j).
    \end{equation}
\end{lemma}

\proof
The proof follows quickly by induction. When $n = 1$, the left-hand side of equality \eqref{eq:combinatorial} only has two terms, \ie
$$\sum_{\vec{\rho} \in \{0,1\}}(-1)^{\rho}(\psi_{\rho}(\xi) - \eta) = (\xi - \eta) - (\bar{\xi} - \eta) =  (\xi - \bar{\xi}),$$
as desired.

Suppose now that equality \eqref{eq:combinatorial} holds for all $n \in \N$ up to some number $N \in \N$. Then, for $n = N + 1$, we deduce using the induction hypothesis that
$$\begin{array}{RCL}
\multicolumn{3}{L}{\sum_{\vec{\rho} \in \{0,1\}^{N + 1}}(-1)^{|\vec{\rho}|}\prod_{j = 1}^{N + 1}(\psi_{\rho_j}(\xi_j) - \eta_j)} \\
~~~~~~ & = & (\xi_{N+1} - \eta_{N+1})\sum_{\substack{\vec{\rho} \in \{0,1\}^{N+1} \\ \rho_{N+1} = 0}}(-1)^{|\vec{\rho}|}\prod_{j = 1}^{N}(\psi_{\rho_j}(\xi_j) - \eta_j) \\
~ & ~ & -\:(\bar{\xi}_{N+1} - \eta_{N+1})\sum_{\substack{\vec{\rho} \in \{0,1\}^{N+1} \\ \rho_{N+1} = 1}}(-1)^{|\vec{\rho}|}\prod_{j = 1}^{N}(\psi_{\rho_j}(\xi_j) - \eta_j) \\
~ & = & \big((\xi_{N+1} - \eta_{N+1}) - (\bar{\xi}_{N+1} - \eta_{N+1})\big)\prod_{j = 1}^{N}(\xi_j - \bar{\xi}_j) = \prod_{j = 1}^{N+1}(\xi_j - \bar{\xi}_j).
\end{array}$$
This finishes the proof.
\endproof

We may now prove our third characterization of the Lebesgue measures, \ie statement (c) of Theorem \ref{thm:Lebesgue_complete}.

\begin{prop}\label{prop:Lebesgue}
    Let $n \in \N$ and let $\mu$ be a Nevanlinna measure in $\R^n$. Then, $\mu = c\lambda_{\R^n}$ for some constant $c \geq 0$ if and only if it satisfies condition \eqref{eq:Lebesgue}, \ie it holds, for every index $j \in \{1,2,\ldots,n\}$, that
    \begin{equation}\tag{\ref{eq:Lebesgue}}
        \int_{\R^n}\frac{1}{(t_j - z_j)^2}\prod_{\substack{\ell = 1 \\ \ell \neq j}}^n\left(\frac{1}{t_\ell - z_\ell} - \frac{1}{t_\ell - \bar{z}_\ell}\right)\diff\mu(\vec{t}) = 0
    \end{equation}
    for all $\vec{z} \in \C^{+n}$.
\end{prop}

\proof
We note first that
$$\int_\R\frac{1}{(t-z)^2}\diff t = 0$$
for every $z \in \C^+$ by standard residue calculus, yielding that condition \eqref{eq:Lebesgue} is satisfied by the measure $c\lambda_{\R^n}$.

Conversely, let $\mu$ be a Nevanlinna measure that satisfies, in addition, condition \eqref{eq:Lebesgue}. Let $q$ be the \HN function given by the data $(0,0,\mu)$ in the sense of Theorem \ref{thm:intRep_Nvar} and let $q_\mathrm{sym}$ be the symmetric extension of the function $q$. Let $j \in \{1,2,\ldots,n\}$ be arbitrary and write
$$f_j(\vec{z}) := \frac{\di}{\di z_j}q_\mathrm{sym}(\vec{z}).$$
We now claim that the function $f_j$ is identically zero on $\C^{+n}$. Without loss of generality, it suffices to show this in the case $j = n$ as all other cases may be considered analogously.

If $n = 1$, it follows, by Proposition \ref{prop:HN_derivatives}, that
$$f_1(z) = \frac{1}{\pi}\int_\R\frac{1}{(t-z)^2}\diff\mu(t) \equiv 0$$
as the measure $\mu$ is assumed to satisfy condition \eqref{eq:Lebesgue}, yielding the desired result. If $n \geq 2$, we proceed as follows. Using Proposition \ref{prop:HN_derivatives} and Lemma \ref{lem:combinatorial}, we deduce, for any $\vec{z} \in \C^{+n}$, that
\begin{multline*}
    \sum_{\vec{\rho} \in \{0,1\}^{n-1}}(-1)^{|\vec{\rho}|}f_n(\psi_{\rho_1}(z_1),\ldots,\psi_{\rho_{n-1}}(z_{n-1}),z_n) \\
    = \frac{2\,\I}{(2\,\pi\,\I)^n}\int_{\R^n}\frac{1}{(t_n-z_n)^2}\bigg[\sum_{\vec{\rho} \in \{0,1\}^{n-1}}(-1)^{|\vec{\rho}|}\prod_{\ell=1}^{n-1}\left(\frac{1}{t_\ell-\psi_{\rho_\ell}(z_\ell)} - \frac{1}{t_\ell + \I}\right)\bigg]\diff\mu(\vec{t}) \\
    = \frac{2\,\I}{(2\,\pi\,\I)^n}\int_{\R^n}\frac{1}{(t_n-z_n)^2}\prod_{\ell=1}^{n-1}\left(\frac{1}{t_\ell-z_\ell} - \frac{1}{t_\ell - \bar{z}_\ell}\right)\diff\mu(\vec{t}) = 0,
\end{multline*}
where the last equality follows from the assumption that the measure $\mu$ satisfies condition \eqref{eq:Lebesgue}. Moreover, we note that, by Proposition \ref{prop:variable_dependence},
$$f_n(\psi_{\rho_1}(z_1),\ldots,\psi_{\rho_{n-1}}(z_{n-1}),z_n) = 0$$
for all index-vectors $\vec{\rho}$ different from $\vec{0}$. As such, it holds that
$$0 = \sum_{\vec{\rho} \in \{0,1\}^{n-1}}(-1)^{|\vec{\rho}|}f_n(\psi_{\rho_1}(z_1),\ldots,\psi_{\rho_{n-1}}(z_{n-1}),z_n) = f_n(z_1,\ldots,z_{n-1},z_n),$$
as desired.

We deduce now that all of the partial derivatives of the function $q_\mathrm{sym}$ are identically zero in $\C^{+n}$, yielding that the same holds for the function $q$. Therefore, the function $q$ is a constant \HN function which, by uniqueness of the representing parameters and Example \ref{ex:Lebesgue}, is only possible if $\mu = c\lambda_{\R^n}$ for some constant $c \geq 0$. This finishes the proof.
\endproof

A slight refinement of the above characterization is the following.

\begin{coro}
    Let $n \in \N$. A Nevanlinna measure $\mu$ on $\R^n$ satisfies condition \eqref{eq:Lebesgue} if and only if it holds, for every index $j \in \{1,2,\ldots,n\}$, that
    $$\int_{\R^n}\frac{1}{(t_j - z_j)^2}\prod_{\substack{\ell = 1 \\ \ell \neq j}}^n\left(\frac{1}{t_\ell - z_\ell} - \frac{1}{t_\ell + \I}\right)\diff\mu(\vec{t}) = 0$$
    for all $\vec{z} \in \C^{+n}$.
\end{coro}

\section{Product Nevanlinna measures}\label{sec:product_measures}

Nevanlinna measures that are formed as a product of the Lebesgue measure and another Nevanlinna measure constitute a very important subclass of Nevanlinna measures.

Firstly, such measures appear as representing measures of \HN functions that do not depend on all of their variables. For example, let $q$ be a \HN function of one variable represented by the data $(a,b,\mu)$ in the sense of Theorem \ref{thm:intRep_Nvar} (for $n = 1$). Then, the \HN function $\til{q}(z_1,z_2) := q(z_1)$ will be represented by the data $(a,(b,0),\mu \otimes \lambda_\R)$ in the sense of Theorem \ref{thm:intRep_Nvar} (for $n = 2$). This follows from the fact that integrating the kernel $K_n$ with respect to \eg $\diff t_n$ gives a constant multiple of the kernel $K_{n-1}$ with the $n$-th variable missing. More precisely, it holds, for every $\vec{z} \in \C^{+n}$ and $\vec{t} \in \R^n$, that
$$\int_\R K_n(\vec{z},\vec{t})\diff t_n = \pi K_{n-1}((z_1,\ldots,z_{n-1}),(t_1,\ldots,t_{n-1})),$$
\cf \cite[Ex. 3.5]{LugerNedic2019} and \cite[Ex. 3.1]{Nedic2017}.

Secondly, product Nevanlinna measures may be interpolated via convex combinations \cite{Nedic2017}. For example, given the measures $\pi\delta_0 \otimes \lambda_\R$ and $\lambda_\R \otimes \pi\delta_0$, where $\delta_0$ denotes the Dirac measure at zero, and a convex combination $k_1 + k_2 = 1$ with $k_1,k_2 > 0$, the measure $\mu$, defined for any Borel set $U \subseteq \R^2$ as
$$\mu(U) := \pi(1+k_2k_1^{-1})\int_\R\chi_U(-k_2k_1^{-1}t,t)\diff t,$$
also is a Nevanlinna measure, \cf \cite[Cor. 4.5]{Nedic2017}, and is the representing measure of the \HN function $(z_1,z_2) \mapsto -(k_1z_1+k_2z_2)^{-1}$. Note that if we allow $k_2 = 0$, the above formula recovers the measure $\lambda_\R \otimes \pi\delta_0$, while for $k_1 = 0$, together with an appropriate change of variables, we recover the measure $\pi\delta_0 \otimes \lambda_\R$. General results on this topic are presented in \eg \cite[Thm. 4.2 and Cor. 4.15]{Nedic2017}.

The following theorem now shows that products with the Lebesgue measure are essentially all product Nevanlinna measures, see also Corollary \ref{coro:mixed_product}.

\begin{thm}\label{thm:product_Nevan}
Let $n_1 + n_2 = n \geq 2$, where $n_1,n_2 \geq 1$. Let $\mu_1$ be a positive Borel measure on $\R^{n_1}$ and $\mu_2$ a positive Borel measure on $\R^{n_2}$. Assume that the measure $\mu := \mu_1 \otimes \mu_2$ is not identically equal to zero. Then, the measure $\mu$ is a Nevanlinna measure if and only if both measures $\mu_1$ and $\mu_2$ are Nevanlinna measures and at least one of them is a constant multiple of the Lebesgue measure (of respective dimension).
\end{thm}

\begin{remark}
If the measure $\mu$ is identically equal to zero, one cannot conclude that both measures $\mu_1$ and $\mu_2$ have to be Nevanlinna measure. Indeed, if $\mu_1 \equiv 0$, then $\mu_2$ may be any positive Borel measure on $\R^{n_2}$.
\end{remark}

\proof
Assume first, without loss of generality, that $\mu_1 = c\lambda_{\R^{n_1}}$ and that $\mu_2$ is a Nevanlinna measure on $\R^{n_2}$. For the the measure $\mu = \mu_1 \otimes \mu_2$, it now follows that 
$$\int_{\R^n}\prod_{j=1}^n\frac{1}{1+t_j^2}\diff\mu(\vec{t}) = \pi^{n_1}\int_{\R^{n_2}}\prod_{j=n_1+1}^{n}\frac{1}{1+t_j^2}\diff\mu_2(t_{n_1+1},\ldots,t_{n}) < \infty,$$
yielding that the measure $\mu$ satisfies the growth condition \eqref{eq:measure_growth}. In order to check that the measure $\mu$ also satisfies the Nevanlinna condition \eqref{eq:measure_Nevan}, we separate three cases with respect to the indices $\ell_1$ and $\ell_2$ in condition \eqref{eq:measure_Nevan}. First, if $\ell_1 \leq n_1$ and $\ell_2 > n_1$, then the integral in condition \eqref{eq:measure_Nevan} is equal to
\begin{multline*}
    \int_\R\frac{1}{(t_{\ell_1}-z_{\ell_1})^2}\diff t_{\ell_1}\cdot(2\,\pi\,\I)^{n_1-1}\\ \cdot  \int_{\R^{n_2}}\frac{1}{(t_{\ell_2} - \bar{z_{\ell_1}})^2}\prod_{\substack{j=n_1+1 \\ j \neq \ell_2}}^{n}\left(\frac{1}{t_j - z_j} - \frac{1}{t_j - \bar{z}_j}\right)\diff\mu_2(t_{n_1+1},\ldots,t_{n}),
\end{multline*}
which is equal to zero for every $\vec{z} \in \C^{+n}$ as
$$\int_\R\frac{1}{(t_{\ell_1}-z_{\ell_1})^2}\diff t_{\ell_1} = 0$$
for any $z_{\ell_1} \in \C^+$ by standard residue calculus. The second case, when $\ell_1,\ell_2 \leq n_1$, follows via an analogous argument. Therefore, it remains to consider the final case, when $\ell_1,\ell_2 > n_1$. In this case, the integral in condition \eqref{eq:measure_Nevan} is equal to
\begin{multline*}
     (2\,\pi\,\I)^{n_1} \\
     \cdot  \int_{\R^{n_2}}\frac{1}{(t_{\ell_1} - z_{\ell_1})^2(t_{\ell_2} - \bar{z_{\ell_1}})^2}\prod_{\substack{j=n_1+1 \\ j \neq \ell_1,\ell_2}}^{n}\left(\frac{1}{t_j - z_j} - \frac{1}{t_j - \bar{z}_j}\right)\diff\mu_2(t_{n_1+1},\ldots,t_{n}).
\end{multline*}
The above integral is precisely the one that appears in the Nevanlinna condition for the measure $\mu_2$ and is, therefore, equal to zero for every $\vec{z} \in \C^{+n}$ by assumption. This finishes the first part of the proof.

Conversely, assume that the measure $\mu = \mu_1 \otimes \mu_2$ is a \NV measure. Then, it satisfies the growth condition \eqref{eq:measure_growth} on $\R^n$, from which we deduce that
$$\int_{\R^{n_1}}\prod_{j=1}^{n_1}\frac{1}{1+t_j^2}\diff\mu_1(t_1,\ldots,t_{n_1})\cdot\int_{\R^{n_2}}\prod_{j=n_1+1}^{n}\frac{1}{1+t_j^2}\diff\mu_2(t_{n_1+1},\ldots,t_{n}) < \infty,$$
implying that both integrals on the left-hand side of the above inequality are finite. Hence, the measures $\mu_1$ and $\mu_2$ satisfy the growth condition \eqref{eq:measure_growth} on $\R^{n_1}$ and $\R^{n_2}$, respectively. 

If $n_1 = 1$, this suffices to conclude that the measure $\mu_1$ is a \NV measure. An analogous reasoning applies if $n_2 = 1$. Otherwise, assume first that $n_1 \geq 2$. Since the measure $\mu$ satisfies the \NV condition \eqref{eq:measure_Nevan} on $\R^n$ and $n_1 \geq 2$, we may investigate what happens when the indices $\ell_1$ and $\ell_2$ in condition \eqref{eq:measure_Nevan} are chosen to satisfy $1 \leq \ell_1 < \ell_2 \leq n_1$ and the vector $\vec{z} \in \C^{+n}$ is chosen such that
$$z_{n_1+1} = \ldots = z_n = \I.$$
If this is the case, it holds that
\begin{multline*}
    \int_{\R^{n_1}}\frac{1}{(t_{\ell_1} - z_{\ell_1})^2(t_{\ell_2} - \bar{z_{\ell_1}})^2}\prod_{\substack{j=1 \\ j \neq \ell_1,\ell_2}}^{n_1}\left(\frac{1}{t_j - z_j} - \frac{1}{t_j - \bar{z}_j}\right)\diff\mu_1(t_1,\ldots,t_{n_1}) \\
    \cdot(2\,\I)^{n_2}\int_{\R^{n_2}}\prod_{j=n_1+1}^n\frac{1}{1+t_j^2}\diff\mu_2(t_{n_1+1},\ldots,t_n) = 0.
\end{multline*}
We note now that the above integral with respect to the measure $\mu_2$ is strictly positive, as the assumption that the measure $\mu$ is not identically zero excludes the case that $\mu_2$ is identically zero. Therefore, the above integral with respect to the measure $\mu_1$ must be zero. From this, we deduce that the measure $\mu_1$ satisfies the \NV condition \eqref{eq:measure_Nevan} on $\R^{n_1}$. If $n_2 \geq 2$, one can show by an analogous procedure that the measure $\mu_2$ satisfies the \NV condition on $\R^{n_2}$. Hence, we conclude that both measures $\mu_1$ and $\mu_2$ are Nevanlinna measures and it remains to show that at least one of them must be equal to a constant multiple of the Lebesgue measure. 

To that end, we consider the Nevanlinna condition \eqref{eq:measure_Nevan} when $1 \leq \ell_1 \leq n_1$ and $n_1 + 1 \leq \ell_2 \leq n$. Taking into account the product structure of the measure $\mu$, we get, for any $\vec{z} \in \C^{+n}$, that
\begin{multline}\label{eq:product_Nevan}
    \int_{\R^{n_1}}\frac{1}{(t_{\ell_1} - z_{\ell_1})^2}\prod_{\substack{j=1 \\ j \neq \ell_1}}^{n_1}\bigg(\frac{1}{t_j - z_j} - \frac{1}{t_j - \bar{z}_j}\bigg)\diff\mu_1(t_1,\ldots,t_{n_1}) \\
    \cdot \int_{\R^{n_2}}\frac{1}{(t_{\ell_2} - \bar{z}_{\ell_1})^2}\prod_{\substack{j=n_1+1 \\ j \neq \ell_2}}^{n}\left(\frac{1}{t_j - z_j} - \frac{1}{t_j - \bar{z}_j}\right)\diff\mu_2(t_{n_1+1},\ldots,t_{n}) = 0.
\end{multline}
We now separate two cases.

\textit{Case 1.} Assume that there exists an index $\ell_2 \in \{n_1 + 1,\ldots,n\}$ and a vector $\vec{\xi} \in \C^{+n_2}$, such that the second integral in equality \eqref{eq:product_Nevan} is non-zero whenever
$$(z_{n_1+1},\ldots,z_n) = \vec{\xi}.$$
If that is the case, it must hold for any index $\ell_1 \in \{1,\ldots,n_1\}$ and any vector $\vec{\zeta} \in \C^{+n_1}$ that
$$\int_{\R^{n_1}}\frac{1}{(t_{\ell_1} - \zeta_{\ell_1})^2}\prod_{\substack{j=1 \\ j \neq \ell_1}}^{n_1}\bigg(\frac{1}{t_j - \zeta_j} - \frac{1}{t_j - \bar{\zeta}_j}\bigg)\diff\mu_1(t_1,\ldots,t_{n_1}) = 0.$$
In other words, the measure $\mu_1$ satisfies condition \eqref{eq:Lebesgue} on $\R^{n_1}$ and is, by Proposition \ref{prop:Lebesgue}, equal to $c\lambda_{\R^{n_1}}$ for some number $c \geq 0$.

\textit{Case 2.} If case 1 does not occur, then it holds, for every index $\ell_2 \in \{n_1 + 1,\ldots,n\}$ and every vector $\vec{\xi} \in \C^{+n_2}$, that
$$\int_{\R^{n_2}}\frac{1}{(t_{\ell_2} - \bar{\xi}_{\ell_1})^2}\prod_{\substack{j=n_1+1 \\ j \neq \ell_2}}^{n}\bigg(\frac{1}{t_j - \xi_j} - \frac{1}{t_j - \bar{\xi}_j}\bigg)\diff\mu_2(t_{n_1+1},\ldots,t_{n}) = 0.$$
Hence, the conjugate of the above equality also holds, leading to an analogous conclusion as in case 1. This finishes the proof.
\endproof

The following corollaries are now an immediate consequence of Theorem \ref{thm:product_Nevan}.

\begin{coro}
If $\sigma$ is a positive Borel measure on $\R^n$ that is not a Nevanlinna measure, there exists no non-zero measure $\mu$ on $\R^m$ such that $\sigma \otimes \mu$ is a Nevanlinna measure on $\R^{n+m}$.
\end{coro}

\begin{coro}
A \HN function whose representing measure is a product measure does non depend on all of its variables.
\end{coro}

Furthermore, Theorem \ref{thm:product_Nevan} may be extended as follows. Let $B := \{k_1,\ldots,k_j\} \subseteq \{1,\ldots,n\}$ be a of size $j := |B| \geq 1$. Define a map $\theta_B\colon \R^n \to \R^j$ as a projection onto the variables indexed by the set $B$, \ie
$$\theta_B(\vec{x}) := (x_{k_1},x_{k_2},\ldots,x_{k_j}).$$
This map has a natural extension $\Theta$ to subsets $V \subseteq \R^n$, defined by
$$\Theta_B(V) := \{\theta_B(\vec{x})~|~\vec{x} \in V\},$$
which allows for the following refinement of the previous result.

\begin{coro}\label{coro:mixed_product}
Let $n_1,n_2$ and $\mu_1,\mu_2$ be as in Theorem \ref{thm:product_Nevan}. Let $B_1 \subseteq \{1,2,\ldots,n\}$ be a set of size $n_1$ and define $B_2 := \{1,2,\ldots,n\} \setminus B_1$. Define a measure $\mu$ on $\R^n$ by setting, for any Borel measurable set $U \subseteq \R^n$,
$$\mu(U) := \mu_1(\Theta_{B_1}(U))\mu_2(\Theta_{B_2}(U)).$$
Then, the measure $\mu$ is a \NV measure on $\R^n$ if and only if both measures $\mu_1$ and $\mu_2$ are Nevanlinna measures and at least one of them is a constant multiple of the Lebesgue measure (of respective dimension).
\end{coro}

\section{The unit polydisk}\label{sec:polydisk}

As the poly-upper half-plane $\C^{+n}$ and the unit polydisk $\D^n$ are biholomorphically equivalent, results may often carry over from one to the other. Therefore, we are first going to recall the definitions of the corresponding classes of functions and measures that relate to the previous sections. Afterwards, the particularities of transitioning to and from the poly-upper half-plane are discussed and, later, the corresponding results concerning characterizations of the Lebesgue measure and product measures are established.

\subsection{Measures with vanishing mixed Fourier coefficients and functions with non-negative real part}\label{sec:polydisk_intro}

The class of measures on $[0,2\pi)^n$, which may also be considered as the unit polytorus $(S^1)^n$, that relates to Nevanlinna measures on $\C^{+n}$ is the following, \cf \cite{Ahern1973,KoranyiPukanszky1963,LugerNedic2018,McDonald1982,McDonald1986,McDonald1990}.

\begin{define}
A finite positive Borel measure on $[0,2\pi)^n$ has \emph{vanishing mixed Fourier coefficients} if
\begin{equation}
    \label{eq:vanishing_mixed_Fourier}
    \int_{[0,2\pi)^n}\E^{\I\,m_1\,s_1}\:\E^{\I\,m_2\,s_2}\,\ldots\,\E^{\I\,m_n\,s_n}\,\diff\nu(\vec{s}) = 0
\end{equation}
for every multi-index $\vec{m} \in \Z^n$ with at least one positive and one negative entry. 
\end{define}

\begin{remark}
If $n = 1$, the multi-index $\vec{m} \in \Z$ only has one entry and thus condition \eqref{eq:vanishing_mixed_Fourier} becomes vacuous. Hence, in this case, one simply considers all finite positive Borel measure on $[0,2\pi)$.
\end{remark}

\begin{remark}
Another class of measures $[0,2\pi)^n$ that has a similar definition is considered in \eg \cite{Forelli1968,Rudin1969}, namely the class of all complex Borel measures on $[0,2\pi)^n$ with vanishing non-negative Fourier coefficents.
\end{remark}

Measures with vanishing mixed Fourier coefficients appear as representing measures for holomorphic functions on $\D^n$ with non-negative real-part via the following theorem, \cf \cite[Thm. 1]{KoranyiPukanszky1963}.

\begin{thm}
    \label{thm:KP}
    A function $f\colon \D^n \to \C$ is holomorphic and has non-negative real part if and only if there exists a number $\alpha \in \R$ and a finite positive Borel measure $\nu$ on $[0,2\pi)^n$ with vanishing mixed Fourier coefficients so that for every $\vec{w} \in \D^n$ it holds that
    \begin{equation}
        \label{eq:intRep_KP}
        f(\vec{w}) = \I\,\alpha + \frac{1}{(2\pi)^n}\int_{[0,2\pi)^n}\bigg(2\prod_{j=1}^n\frac{1}{1-w_j\,\E^{-\I\,s_j}}-1\bigg)\diff\nu(\vec{s}).
    \end{equation}
    Furthermore, the correspondence between the pair ($\alpha$,$\nu$) and the function $f$ is unique.
\end{thm}

\begin{example}\label{ex:Lebesgue_polydisk}
Let $n \in \N$ and consider the function $f\colon \D^n \to \C$ given by
$$f(\vec{w}) := 1.$$
This function is represented by the data $(0,\lambda_{[0,2\pi)^n})$ in the sense of representation \eqref{eq:intRep_KP}, where $\lambda_{[0,2\pi)^n}$ denotes the standard Lebesgue measure on $[0,2\pi)^n$. More generally, the function
$$f(\vec{w}) := \xi$$
with $\Re[\xi] \geq 0$ is represented by the data $(\Im[\xi],\Re[\xi]\lambda_{[0,2\pi)^n})$.\hfill$\lozenge$
\end{example}

\subsection{Transitioning to and from the poly-upper half-plane}\label{sec:transition}

Using the Cayley transform and its inverse, any holomorphic function on the unit polydisk may be transformed to a function on poly-upper half-plane. We recall that the Cayley transform $\varphi\colon \C^+ \to \D$ and its inverse $\varphi^{-1}\colon \D \to \C^+$ are defined as
$$\varphi(z) := \frac{z-\I}{z+\I} \quad\text{and}\quad \varphi^{-1}(w) := \I\,\frac{1+w}{1-w},$$
where $z \in \C^+$ and $w \in \D$. As the Cayley transform and its inverse also constitute bijections between $\R$ and $S^1\setminus\{1\}$, we may use them to change variables when integrating \eg over $[0,2\pi)$ by setting $\E^{\I\,s} = \varphi(t)$ for $t \in \R$ and $s \in [0,2\pi)$. This change of variables also contributes a Jacobian factor of $2(1+t^2)^{-1}$ or $(1-\cos(s))^{-1}$ depending on whether we are changing the $s$-variable to the $t$-variable or the other way around. We note also that some authors may define the map $\varphi^{-1}$ as the Cayley transform instead.

If we now have a \HN function $q$, we may assign to it, via the inverse Cayley transform, a function $f$ on $\D^n$ with non-negative real part by setting
\begin{equation}
    \label{eq:q_to_f}
    f(\vec{w}) := -\I^{-1}q(\varphi^{-1}(w_1),\varphi^{-1}(w_2),\ldots,\varphi^{-1}(w_n)).
\end{equation}
Conversely, starting with such a function $f$, we arrive at a \HN function $q$ by setting
\begin{equation}
    \label{eq:f_to_q}
    q(\vec{z}) := \I\,f(\varphi(z_1),\varphi(z_2),\ldots,\varphi(z_n)).
\end{equation}
Doing this process twice always gives gives back the original function. Therefore, when we say that a function $f$ on $\D^n$ corresponds to a function $q$ on $\C^{+n}$, or vice-versa, this is meant in the sense that one is given in terms of the other by one of the formulas \eqref{eq:q_to_f} and \eqref{eq:f_to_q}. Likewise, we say that a data-set $(a,b,\mu)$ from Theorem \ref{thm:intRep_Nvar} corresponds to a data-pair $(\alpha,\nu)$ from Theorem \ref{thm:KP} if the functions they represent correspond to one-another in the above meaning.

\begin{example}
We may infer immediately from relations \eqref{eq:q_to_f} and \eqref{eq:f_to_q} that the functions $q$ and $f$ from Examples \ref{ex:Lebesgue} and \ref{ex:Lebesgue_polydisk} correspond to one-another. Hence, the data-set $(0,\vec{0},\lambda_{\R^n})$ corresponds to the data-pair $(0,\lambda_{[0,2\pi)^n})$. In particular, we may say that the Lebesgue measure on $\R^n$ is corresponds to the Lebesgue measure on $[0,2\pi)^n$.\hfill$\lozenge$
\end{example}

It is important to note that, in general, changing the measure $\nu$ in a pair $(\alpha,\nu)$ may change both the measure $\mu$ and the vector $\vec{b}$ in the corresponding set $(a,\vec{b},\mu)$. Conversely, changing the vector $\vec{b}$ in a data-set $(a,\vec{b},\mu)$ will change the measure $\nu$ in the corresponding pair $(\alpha,\nu)$. This follows from fact that the restriction of a measure $\nu$ with vanishing mixed Fourier coefficients to a set of the form $\{\vec{s} \in [0,2\pi)^n~|~s_j = 0\}$ for some $j \in \{1,\ldots,n\}$ is equal to a constant multiple of $\lambda_{[0,2\pi)^{n-1}}$, \cf \cite[Cor. 3.7]{LugerNedic2017} and \cite[Thm. 4.1]{LugerNedic2018}, see also \cite[Thm. 4.1]{LugerNedic2019}. However, we note that changing the number $\alpha$ only affects the number $a$ and vice-versa.

\subsection{Characterizations of the Lebesgue measure}\label{sec:polydisk_characterizations}

A characterization of the Lebesgue measure on $[0,2\pi)^n$ in terms of its Fourier coefficients is well known and follows from the injectivity of the Fourier transform, see Remark \ref{rem:Fourier} later on. However, the result may also be obtained using only representation \eqref{eq:intRep_KP} and we provide the proof for the sake of completeness.

\begin{prop}
    \label{prop:Lebesgue_polydisk}
    Let $\nu$ be a finite positive Borel measure on $[0,2\pi)^n$ with vanishing mixed Fourier coefficients. Then, $\nu = c\,\lambda_{[0,2\pi)^n}$ for some constant $c \geq 0$ if and only if the measure $\nu$ also has vanishing non-zero non-negative Fourier coefficients, \ie 
    \begin{equation}
    \label{eq:vanishing_non_zero_Fourier}
    \int_{[0,2\pi)^n}\E^{\I\,m_1\,s_1}\:\E^{\I\,m_2\,s_2}\,\ldots\,\E^{\I\,m_n\,s_n}\,\diff\nu(\vec{s}) = 0
    \end{equation}
    for every multi-index $\vec{m} \in (\N_0)^n\setminus\{\vec{0}\}$.
\end{prop}

\proof
Note first that
$$\int_{[0,2\pi)} \E^{\I\,m\,s}\diff s = 0$$
for every $m \neq 0$, yielding that condition \eqref{eq:vanishing_non_zero_Fourier} is satisfied by the measure $c\lambda_{[0,2\pi)^n}$.

Conversely, let $\nu$ be a measure with vanishing mixed Fourier coefficients that satisfies, in addition, condition \eqref{eq:vanishing_non_zero_Fourier} and let $f$ be the function on $\D^n$ given by the pair $(0,\nu)$ in the sense of Theorem \ref{thm:KP}.

We are now going to show that all partial derivatives of the function $f$ are identically zero. Without loss of generality, it suffices to show this only for the derivative with respect to the variable $w_1$, as all other cases follow an analogous argument. To that end, we calculate that
\begin{equation}
    \label{eq:polydisk_derivative}
    \frac{\di f}{\di w_1}(\vec{w}) = \frac{1}{(2\pi)^n}\int_{[0,2\pi)^n}\frac{2\,\E^{-\I\,s_1}}{(1-w_1\,\E^{-\I\,s_1})^2}\prod_{j=2}^n\frac{1}{1-w_j\,\E^{-\I\,s_j}}\diff\nu(\vec{s}),
\end{equation}
where we are allowed to move the derivative under the integral sign as the integrand is bounded function in the $\vec{s}$-variables whenever $\vec{w}$ is restricted to a compact subset of $\D^n$.

If $\vec{w}$ remains restricted to a compact subset of $\D^n$, we may use a geometric series expansion to rewrite the integrand in right-hand side of equality \eqref{eq:polydisk_derivative}, yielding that
\begin{eqnarray*}
\frac{\di f}{\di w_1}(\vec{w}) & = & \frac{2}{(2\pi)^n}\int_{[0,2\pi)^n}\bigg(\sum_{\substack{\vec{m} \in (\N_0)^n \\ m_1 \neq 0}}m_1\,\E^{-\I\,m_1\,s_1}\ldots\E^{-\I\,m_n\,s_n}w_1^{m_1}\ldots w_n^{m_n}\bigg)\diff\nu(\vec{s}), \\
~ & = & \frac{2}{(2\pi)^n}\sum_{\substack{\vec{m} \in (\N_0)^n \\ m_1 \neq 0}}m_1\,w_1^{m_1}\ldots w_n^{m_n}\int_{[0,2\pi)^n}\E^{-\I\,m_1\,s_1}\,\ldots\,\E^{-\I\,m_n\,s_n}\,\diff\nu(\vec{s}) \\
~ & = & 0.
\end{eqnarray*}
As the compact set restricting $\vec{w}$ can be chosen arbitrary, we conclude that
$$\frac{\di f}{\di w_1}(\vec{w}) = 0$$
for every $\vec{w} \in \D^n$. This finishes the proof.
\endproof

Proposition \ref{prop:Lebesgue_polydisk} gives, after transforming condition \eqref{eq:vanishing_non_zero_Fourier} via the Cayley transform, the last part of the proof of Theorem \ref{thm:Lebesgue_complete}, as it implies that condition (d) also characterizes the Lebesgue measure and is, hence, equivalent to the other conditions in the theorem. 

Similarly, using the inverse Cayley transform, we may translate the statement of Theorem \ref{thm:Lebesgue_complete} to yield the three corresponding characterizations of the Lebesgue measure among all measures with vanishing mixed Fourier coefficients. This requires, also, that we translate the factors $N_{j,k}$ by setting $t = \varphi^{-1}(\E^{\I\,s})$ and $z = \varphi^{-1}(w)$. By also counting in the Jacobian factor $(1-\cos(s))^{-1}$, we define, given ambient vectors $\vec{s} \in [0,2\pi)^n$ and $\vec{w} \in \D^n$, the factors $D_{k,j}$ for $k \in \{-1,0,1\}$ and $j \in \{1,2,\ldots,n\}$ as
$$D_{-1,j} := \frac{1}{2}\,\frac{w_j}{\E^{\I\,s_j}-w_j}, \quad D_{0,j} := \frac{1}{2}, \quad
D_{1,j} := -\frac{1}{2}\,\frac{\bar{w}_j}{\E^{-\I\,s_j}-\bar{w}_j}.$$

The corresponding analogue of Theorem \ref{thm:Lebesgue_complete} for measures with vanishing mixed Fourier coefficients is, therefore, the following.

\begin{thm}\label{thm:Lebesgue_complete_polydisk}
Let $\nu$ be a positive Borel measure on $[0,2\pi)^n$ with vanishing mixed Fourier coefficients. Then, $\nu = c\lambda_{[0,2\pi)^n}$ for some constant $c \geq 0$ if and only if the measure $\nu$ satisfies one (and thus all) of the following conditions:
\begin{itemize}
    \item[(a)]{It holds that
    $$\sum_{\substack{\vec{\rho} \in \{-1,0,1\}^n \\ 1\not\in\vec{\rho}\,\wedge\,\vec{\rho} \neq \vec{0}}}\int_{[0,2\pi)^n}D_{\rho_1,1}D_{\rho_2,2}\ldots D_{\rho_n,n}\diff\nu(\vec{s}) = 0$$
    for all $\vec{w} \in \D^n$.    
    }
    \item[(b)]{It holds that
    $$\int_{[0,2\pi)^n}D_{\rho_1,1}D_{\rho_2,2}\ldots D_{\rho_n,n}\diff\nu(\vec{s}) = 0$$
    for all $\vec{w} \in \D^n$ and and for every vector $\vec{\rho} \in \{-1,0,1\}^n$ that is not the zero-vector and has no entries equal to $1$.
    }
    \item[(c)]{The measure $\nu$ satisfies condition \eqref{eq:vanishing_non_zero_Fourier}, \ie it holds that
    $$\int_{[0,2\pi)^n}\E^{\I\,m_1\,s_1}\:\E^{\I\,m_2\,s_2}\,\ldots\,\E^{\I\,m_n\,s_n}\,\diff\nu(\vec{s}) = 0$$
    for every multi-index $\vec{m} \in (\N_0)^n\setminus\{\vec{0}\}$.
    } 
    \item[(d)]{It holds, for every index $j \in \{1,2,\ldots,n\}$, that
    $$\int_{[0,2\pi)^n}\E^{\I\,s_j}\left(\frac{w_j-1}{w_j - \E^{\I\,s_j}}\right)^2\prod_{\substack{\ell = 1 \\ \ell \neq j}}^{n}\left(\frac{w_\ell}{\E^{\I\,s_\ell} - w_\ell} - \frac{\E^{-\I\,s_\ell}}{\E^{-\I\,s_\ell}-\bar{w}_\ell}\right)\diff\nu(\vec{s}) = 0$$
    for all $\vec{w} \in \D^n$.
    } 
\end{itemize}
\end{thm}

A reformulation of the above theorem may be stated as follows

\begin{coro}\label{coro:Lebesgue_polydisk}
Let $\nu$ be a positive Borel measure on $[0,2\pi)^n$. Then, $\nu = c\lambda_{[0,2\pi)^n}$ for some constant $c \geq 0$ if and only if the measure $\nu$ satisfies one (and thus all) of the following conditions:
\begin{itemize}
    \item[(a)]{It holds that
    $$\sum_{\substack{\vec{\rho} \in \{-1,0,1\}^n \\ \vec{\rho} \neq \vec{0}}}\int_{[0,2\pi)^n}D_{\rho_1,1}D_{\rho_2,2}\ldots D_{\rho_n,n}\diff\nu(\vec{s}) = 0$$
    for all $\vec{w} \in \D^n$.    
    }
    \item[(b)]{It holds that
    $$\int_{[0,2\pi)^n}D_{\rho_1,1}D_{\rho_2,2}\ldots D_{\rho_n,n}\diff\nu(\vec{s}) = 0$$
    for all $\vec{w} \in \D^n$ and and for every vector $\vec{\rho} \in \{-1,0,1\}^n$ that is not the zero-vector.
    }
    \item[(c)]{All of the non-zero Fourier coefficients of $\nu$ are zero \ie it holds that
    $$\int_{[0,2\pi)^n}\E^{\I\,m_1\,s_1}\:\E^{\I\,m_2\,s_2}\,\ldots\,\E^{\I\,m_n\,s_n}\,\diff\nu(\vec{s}) = 0$$
    for every multi-index $\vec{m} \in \Z^n\setminus\{\vec{0}\}$.
    } 
    \item[(d)]{It holds, for every index $j \in \{1,2,\ldots,n\}$, that
    $$\int_{[0,2\pi)^n}\E^{\I\,s_j}\left(\frac{w_j-1}{w_j - \E^{\I\,s_j}}\right)^2\prod_{\substack{\ell = 1 \\ \ell \neq j}}^{n}\left(\frac{w_\ell}{\E^{\I\,s_\ell} - w_\ell} - \frac{\E^{-\I\,s_\ell}}{\E^{-\I\,s_\ell}-\bar{w}_\ell}\right)\diff\nu(\vec{s}) = 0$$
    for every $\vec{w} \in \D^n$ and it holds, for all indices $j_1,j_2 \in \{1,2,.\ldots,n\}$ with $j_1 < j_2$, that
    \begin{multline*}
        \int_{[0,2\pi)^n}\E^{\I\,s_{j_1}}\left(\frac{w_{j_1}-1}{w_{j_1} - \E^{\I\,s_{j_1}}}\right)^2\E^{\I\,s_{j_2}}\left(\frac{\bar{w}_{j_2}-1}{\bar{w}_{j_2}\E^{\I\,s_{j_2}} - 1}\right)^2 \\ \cdot\prod_{\substack{\ell = 1 \\ \ell \neq j_1,j_2}}^{n}\left(\frac{w_\ell}{\E^{\I\,s_\ell} - w_\ell} - \frac{\E^{-\I\,s_\ell}}{\E^{-\I\,s_\ell}-\bar{w}_\ell}\right)\diff\nu(\vec{s}) = 0
    \end{multline*}
    for every $\vec{w} \in \D^n$.
    } 
\end{itemize}
\end{coro}

\begin{remark}\label{rem:Fourier}
Due to the injectivity of the Fourier transform, see \eg \cite[Thm. 2.1]{SteinShakarchi2003} for the classical setting or \cite[Thm. 1.9.5]{Sasvari1994} for the abstract setting, we know that if all of the Fourier coefficients of a (complex) Borel measure $\nu$ on $[0,2\pi)^n$ are zero, then $\nu \equiv 0$. If, instead, all but the coefficient with multi-index $\vec{0}$ are zero, we may consider the measure $\til{\nu} := \nu - c\lambda_{[0,2\pi)^n}$. Then, it holds that all of the Fourier coefficients of $\til{\nu}$ are zero, implying that $\til{\nu} \equiv 0$ and yielding the same characterization as Proposition \ref{prop:Lebesgue_polydisk} or statement (c) of Corollary \ref{coro:Lebesgue_polydisk}.
\end{remark}

\subsection{Product measures with vanishing mixed Fourier coefficients}

An analogous result concerning product measures with vanishing mixed Fourier coefficients may now also be established.

\begin{prop}
    \label{prop:product_polydisk}
    Let $n_1 + n_2 = n \geq 2$, where $n_1,n_2 \geq 1$. Let $\nu_1$ be a finite positive Borel measure on $[0,2\pi)^{n_1}$ and $\nu_2$ a positive Borel measure on $[0,2\pi)^{n_2}$. Assume that the measure $\nu := \nu_1 \otimes \nu_2$ is not identically equal to zero. Then, the measure $\nu$ has vanishing mixed Fourier coefficients if and only if both measures $\nu_1$ and $\nu_2$ have vanishing mixed Fourier coefficients and at least one of them is a constant multiple of the Lebesgue measure (of respective dimension).
\end{prop}

\proof
Assume first, without loss of generality, that $\nu_1 = d\lambda_{[0,2\pi)^{n_1}}$ and that $\nu_2$ is a finite positive Borel measure on $\R^{n_2}$ with vanishing mixed Fourier coefficients. Then, the measure $\nu := \nu_1 \otimes \nu_2$ is finite positive Borel measure on $\R^n$ and its Fourier coefficients are given by the product
\begin{multline}
    \label{eq:poly_product}
    \int_{[0,2\pi)^{n_1}}\E^{\I\,m_1\,s_1}\ldots\E^{\I\,m_{n_1}\,s_n}\diff s_1,\ldots \diff s_{n_1} \\
    \cdot \int_{[0,2\pi)^{n_2}}\E^{\I\,m_{n_1+1}\,s_{n_1+1}}\ldots\E^{\I\,m_n\,s_n}\diff\nu_2(s_{n_1+1},\ldots s_n).
\end{multline}
Take now any multi-index $\vec{m} \in \Z^n$ that has at least one positive and at least one negative entry. If there is either a positive or negative entry among the positions 1 to $n_1$, then the corresponding Fourier coefficient is zero as the first factor in the product \eqref{eq:poly_product} is zero. Otherwise, if all of the first $n_1$ entries of $\vec{m}$ are zero, there must be at least one positive and at least one negative entry appearing among the positions $n_1 + 1$ to $n$. In this case, the second factor in the product \eqref{eq:poly_product} is a mixed Fourier coefficient of the measure $\nu_2$, which is zero by assumption. Thus, the measure $\nu$ has vanishing mixed Fourier coefficients.

Conversely, assume that the measure $\nu := \nu_1 \otimes \nu_2$ is a finite positive Borel measure on $[0,2\pi)^{n}$ with vanishing mixed Fourier coefficients. Then, the measures $\nu_1$ and $\nu_2$ must also be finite positive Borel measures as the measure $\nu$ is assumed to not be identically zero. To see that the measure $\nu_1$ must have vanishing mixed Fourier coefficients, let the multi-index $\vec{m} \in \Z^n$ be such that it has at least one positive and at least one negative entry among the first $n_1$ entries with
$$m_{n_1+1} = \ldots = m_n = 0.$$
Then, the Fourier coefficient of the measure $\nu$ with index $\vec{m}$ is equal to zero by assumption, implying that
$$\int_{[0,2\pi)^{n_1}}\E^{\I\,m_1\,s_1}\ldots\E^{\I\,m_{n_1}\,s_n}\diff\nu_1(s_1,\ldots,s_{n_1}) \cdot \int_{[0,2\pi)^{n_2}}\diff\nu_2(s_{n_1+1},\ldots,s_n) = 0.$$
As the measure $\nu_2$ is not identically zero, we deduce that the measure $\nu_1$ has vanishing mixed Fourier coefficients. An analogous reasoning may be done for the measure $\nu_2$.

Hence, it remains to show that at least one of the measures $\nu_1$ and $\nu_2$ is equal to a constant multiple of the Lebesgue measure, \ie at least one of them satisfies condition \eqref{eq:vanishing_non_zero_Fourier}. To that end, we consider two cases.

\textit{Case 1.} Assume that there exists a multi-index $\vec{m} \in \Z^n$ such that all entries at positions $n_1 + 1$ to $n$ are non-positive with at least one being strictly negative, such that
$$\int_{[0,2\pi)^{n_2}}\E^{\I\,m_{n_1+1}\,s_{n_1+1}}\ldots\E^{\I\,m_n\,s_n}\diff\nu_2(s_{n_1+1},\ldots s_n) \neq 0.$$
Then, the Fourier coefficient of the measure $\nu$ is zero as soon as the multi-index $\vec{m} \in \Z^n$ has, in addition, only non-negative entries among the first $n_1$ entries with at least one entry being strictly positive, yielding that
$$\int_{[0,2\pi)^{n_1}}\E^{\I\,m_1\,s_1}\ldots\E^{\I\,m_{n_1}\,s_n}\diff\nu_1(s_1,\ldots,s_{n_1}) = 0$$
for any such multi-index $\vec{m}$. The measure $\nu_1$ satisfies, therefore, condition \eqref{eq:vanishing_non_zero_Fourier}, as desired.

\textit{Case 2.} If case 1 does not occur, then it must hold, for every multi-index $\vec{m} \in \Z^n$ such that all entries at positions $n_1 + 1$ to $n$ are non-positive with at least one being strictly negative that
$$\int_{[0,2\pi)^{n_2}}\E^{\I\,m_{n_1+1}\,s_{n_1+1}}\ldots\E^{\I\,m_n\,s_n}\diff\nu_2(s_{n_1+1},\ldots s_n) = 0.$$
Hence, the conjugate of the above equality also holds, implying that the measure $\nu_2$ satisfies condition condition \eqref{eq:vanishing_non_zero_Fourier}, finishing the proof.
\endproof

Note that Proposition \ref{prop:product_polydisk} does not follow directly from Theorem \ref{thm:product_Nevan}. Indeed, if $f$ is a holomorphic function on $\D^n$ represented by the pair $(0,\nu_1 \otimes \nu_2)$ in the sense of representation \eqref{eq:intRep_KP}, then the function $q$ constructed via formula \eqref{eq:f_to_q} will be a \HN function, represented by some data $(0,\vec{b},\mu)$ in the sense of representation \eqref{eq:intRep_Nvar}. However, there is \emph{a priori} no reason to assume that this measure $\mu$ is equal to the measure $\mu_1 \otimes \mu_2$, obtained by transforming the measure $\nu_1 \otimes \nu_2$ via the Cayley transform as usual, though this may be inferred from \eg the proof of \cite[Thm. 4.1]{LugerNedic2019}. Therefore, the proof presented above is preferred as it is self-contained.

Analogous corollaries to those of Theorem \ref{thm:product_Nevan} may now also be presented.

\begin{coro}
If $\sigma$ is a finite positive Borel measure on $[0,2\pi)^n$ that does not have vanishing mixed Fourier coefficients, there exists no non-zero measure $\nu$ on $[0,2\pi)^n$ such that the measure $\sigma \otimes \nu$ has vanishing mixed Fourier coefficients.
\end{coro}

\begin{coro}
A holomorphic function on $\D^n$ with non-negative real part whose representing measure is a product measure does not depend on all of its variables.
\end{coro}

\begin{coro}
Let $n_1,n_2$ and $\nu_1,\nu_2$ be as in Proposition \ref{prop:product_polydisk}. Let $B_1 \subseteq \{1,2,\ldots,n\}$ be a set of size $n_1$ and define $B_2 := \{1,2,\ldots,n\} \setminus B_1$. Define a measure $\nu$ on $[0,2\pi)^n$ via the map $\Theta$ from Corollary \ref{coro:mixed_product}, \ie by setting, for any Borel measurable set $U \subseteq [0,2\pi)^n$,
$$\nu(U) := \nu_1(\Theta_{B_1}(U))\nu_2(\Theta_{B_2}(U)).$$
Then, the measure $\nu$ has vanishing mixed Fourier coefficients if and only if both measures $\nu_1$ and $\nu_2$ have vanishing mixed Fourier coefficients and at least one of them is a constant multiple of the Lebesgue measure (of respective dimension).
\end{coro}

\section*{Acknowledgments}

The author would like to thank Annemarie Luger for enthusiastic discussions on the topic and Alan Sola for careful reading of the manuscript and many helpful comments and suggestions.

\bibliographystyle{amsplain}
\bibliography{total}

\end{document}